\documentclass[11pt]{amsart}
\usepackage{times,epsf,amssymb,amsmath,hyperref}

\begin{document}

\newtheorem{thm}{Theorem}[section]
\newtheorem{lem}[thm]{Lemma}
\newtheorem{cor}[thm]{Corollary}

\theoremstyle{definition}
\newtheorem{defn}{Definition}[section]

\theoremstyle{remark}
\newtheorem{rmk}{Remark}[section]

\def\square{\hfill${\vcenter{\vbox{\hrule height.4pt \hbox{\vrule
width.4pt height7pt \kern7pt \vrule width.4pt} \hrule height.4pt}}}$}

\def\T{\mathcal T}

\newenvironment{pf}{{\it Proof:}\quad}{\square \vskip 12pt}

\title{Generic Uniqueness of Area Minimizing Disks for Extreme Curves}
\author{Baris Coskunuzer}
\address{Koc University \\ Department of Mathematics \\ Sariyer, Istanbul 34450 Turkey}
\email{bcoskunuzer@ku.edu.tr}
\thanks{The author is partially supported by EU-FP7 Grant IRG-226062 and TUBITAK Grant 107T642}

\maketitle


\newcommand{\SH}{S^2_{\infty}(\mathbf{H}^3)}
\newcommand{\PM}{\partial M}
\newcommand{\SI}{S^2_{\infty}}
\newcommand{\BHH}{\mathbf{H}^3}
\newcommand{\CH}{\mathcal{C}(\Gamma)}
\newcommand{\BH}{\mathbf{H}}
\newcommand{\BR}{\mathbf{R}}
\newcommand{\BC}{\mathbf{C}}
\newcommand{\BZ}{\mathbf{Z}}

\begin{abstract}

We show that for a generic nullhomotopic simple closed curve $\Gamma$ in the boundary of a compact, orientable, mean
convex $3$-manifold $M$ with $H_2(M,\BZ)=0$, there is a unique area minimizing disk $D$ embedded in $M$ with $\partial
D = \Gamma$. We also show that the same is true for nullhomologous curves in absolutely area minimizing surface case.

\end{abstract}

\section{Introduction}

The Plateau problem asks the existence of an area minimizing disk for a given curve in the ambient manifold $M$. This
problem was solved for $\BR^3$ by Douglas \cite{Do}, and Rado \cite{Ra} in early 1930s. Later, it was generalized by
Morrey \cite{Mo} for Riemannian manifolds. Then, regularity (nonexistence of branch points) of these solutions was
shown by Osserman \cite{Os}, Gulliver \cite{Gu} and Alt \cite{Al}. In the early 1960s, the same question was studied
for absolutely area minimizing surfaces, i.e. for surfaces that minimize area among all oriented surfaces with the
given boundary (without restriction on genus). The geometric measure theory techniques proved to be quite powerful, and
De Georgi, Federer-Fleming solved the problem for area minimizing surfaces \cite{Fe}.

Later, the question of embeddedness of the solution was studied by many experts. First, Gulliver-Spruck showed
embeddedness for the extreme curves with total curvature less than $4\pi$ in \cite{GS}. Tomi-Tromba \cite{TT} and
Almgren-Simon \cite{AS} showed the existence of embedded minimal (not necessarily area minimizing) disks for extreme
curves. Then, Meeks-Yau \cite{MY1} showed that, for extreme boundary curves, area minimizing disks must be embedded.
Recently, Ekholm, White, and Wienholtz generalized Gulliver-Spruck embeddedness result by removing extremeness
condition from the curves \cite{EWW}.

On the other hand, the number of the solutions was also an active area of research. First, Rado showed that if a curve
can be projected bijectively to a convex plane curve, then it bounds a unique minimal disk. Then, Nitsche proved
uniqueness of minimal disks for the boundary curves with total curvature less than $4\pi$ in \cite{Ni}. Then, Tromba
\cite{Tr} showed that a generic curve in $\BR^3$ bounds a unique area minimizing disk. Then, Morgan \cite{M} proved a
similar result for area minimizing surfaces. Later, White proved a very strong generic uniqueness result for fixed
topological type in any dimension \cite{Wh1}. In particular, he showed that a generic $k$-dimensional, $C^{j,\alpha}$
submanifold of a Riemannian manifold cannot bound two smooth, minimal $(k+1)$-manifolds having the same area.

In this paper, we will give a new generic uniqueness results for both versions of the Plateau problem. Our techniques
are simple and topological. The first main result is the following:

\vspace{.1in}

\noindent {\bf Theorem 3.2:} Let $M$ be a compact, orientable, mean convex $3$-manifold with $H_2(M,\BZ)=0$. Then for a
generic nullhomotopic (in $M$) simple closed curve $\Gamma$ in $\partial M$, there exists a unique area minimizing disk
$D$ in $M$ with $\partial D = \Gamma$.

\vspace{.1in}

This theorem is also true for compact, irreducible, orientable, mean convex $3$-manifolds (See Remark 3.2). The second
main result is a similar theorem for absolutely area minimizing surfaces.

\vspace{.1in}

\noindent {\bf Theorem 4.3:} Let $M$ be a compact, orientable, mean convex $3$-manifold with $H_2(M,\BZ)=0$. Then for a
generic nullhomologous (in $M$) simple closed curve $\Gamma$ in $\partial M$, there exists a unique absolutely area
minimizing surface $\Sigma$ in $M$ with $\partial \Sigma = \Gamma$.

\vspace{.1in}

These results naturally generalize to noncompact homogeneously regular $3$-manifolds (see the last section).

The short outline of the technique for generic uniqueness is the following: For simplicity, we will
focus on the case of the area minimizing disks in a mean convex manifold $M$. Let $\Gamma_0$ be a
nullhomotopic (in $M$) simple closed curve in $\PM$. First, we will show that either there exists a
unique area minimizing disk $ D_0$ in $M$ with $\partial D_0=\Gamma_0$, or there exist two {\em
disjoint} area minimizing disks $ D_0^+ ,  D_0^-$ in $M$ with $\partial D_0^\pm=\Gamma_0$.

Now, take a small neighborhood $N(\Gamma_0)\subset \PM$ which is an annulus. Then foliate
$N(\Gamma_0)$ by simple closed curves $\{\Gamma_t\}$ where $t\in(-\epsilon, \epsilon)$, i.e.
$N(\Gamma_0) \simeq \Gamma\times (-\epsilon, \epsilon)$. By the above fact, for any $\Gamma_t$
either there exists a unique area minimizing disk $ D_t$, or there are two area minimizing disks $
D_t^\pm$ disjoint from each other. Also, since these are area minimizing disks, if they have
disjoint boundary, then they are disjoint by \cite{MY2}. This means, if $t_1<t_2$, then $ D_{t_1}$
is disjoint and {\em below} $ D_{t_2}$ in $M$. Consider this collection of area minimizing disks.
Note that for curves $\Gamma_t$ bounding more than one area minimizing disk, we have a canonical
region $N_t$ in $M$ between the disjoint area minimizing disks $ D_t^\pm$.

Now, take a finite curve $\beta\subset M$ which is transverse to the collection of these area
minimizing disks $\{ D_t\}$ whose boundaries are $\{\Gamma_t\}$. Let the length of this line
segment be $C$.

Now, the idea is to consider the {\em thickness} of the neighborhoods $N_t$ assigned to the
boundary curves $\{\Gamma_t\}$. When $\Gamma_t$ bounds a unique area minimizing disk $D_t$, let
$N_t=D_t$ be a degenerate canonical region for $\Gamma_t$. Let $s_t$ be the length of the segment
$I_t$ of $\beta$ between $ D_t^+$ and $ D_t^-$, which is the {\em width} of $N_t$ assigned to
$\Gamma_t$. Then, the curves $\Gamma_t$ bounding more than one area minimizing disk have positive
width, and contributes to total thickness of the collection, and the curves bounding a unique area
minimizing disk has $0$ width and do not contribute to the total thickness. Since
$\sum_{t\in(-\epsilon, \epsilon)} s_t < C$, the total thickness is finite. This implies for only
countably many $t\in(-\epsilon, \epsilon)$, $s_t>0$, i.e. $\Gamma_t$ bounds more than one area
minimizing disk. For the remaining uncountably many $t\in(-\epsilon, \epsilon)$, $s_t=0$, and there
exists a unique area minimizing disk for those $t$. This proves the space of simple closed curves
of uniqueness is dense in the space of Jordan curves in $\PM$. Then, we will show this space is not
only dense, but also generic.

\vspace{.1in}

The organization of the paper is as follows: In the next section we will cover some basic results
which will be used in the following sections. In section 3, we will prove the first main result of
the paper. Then in section 4, we will show the area minimizing surfaces case. Finally in section 5,
we will have some final remarks.

\subsection{Acknowledgements:}

I am very grateful to the referee for very valuable comments and suggestions. I would like to thank Brian White and
Frank Morgan for very useful conversations.

\section{Preliminaries}

In this section, we will overview the basic results which we use in the following sections. First, we should note that
Hass-Scott's very nicely written paper \cite{HS} would be a great reference for a good introduction for the notions in
this paper. We will start with some basic definitions.

\begin{defn} An {\em area minimizing disk} is a disk which has the smallest area among the disks with the same boundary.
An {\em absolutely area minimizing surface} is a surface which has the smallest area among all orientable surfaces
(with no topological restriction) with the same boundary.
\end{defn}

\begin{defn} Let $M$ be a compact Riemannian $3$-manifold with boundary. Then $M$ is a {\em mean convex} (or sufficiently convex) if the following conditions hold.

\begin{itemize}

\item $\partial M$ is piecewise smooth.

\item Each smooth subsurface of $\partial M$ has nonnegative curvature with respect to inward normal.

\item There exists a Riemannian manifold $N$ such that $M$ is isometric to a submanifold of $N$ and
each smooth subsurface $S$ of $\partial M$  extends to a smooth embedded surface $S'$ in $N$ such
that $S' \cap M = S$.

\end{itemize}

\end{defn}

\begin{defn} A simple closed curve is an {\em extreme curve} if it is on the boundary of its convex hull.
A simple closed curve is called as {\em $H$-extreme curve} if it is a curve in the boundary of a mean convex manifold $M$.
\end{defn}

\begin{rmk}
Note that our results in this paper are for $H$-extreme curves which are in the boundary of a fixed $3$-manifold $M$.
Since any extreme curve is also $H$-extreme, our results applies to this case as well. Note also that for any smooth
embedded curve $\Gamma$, one can find a mean convex (sufficiently thin) solid torus $T_\Gamma$ such that $\Gamma\subset
\partial T_\Gamma$, hence $\Gamma$ is $H$-extreme. So, an $H$-extreme curve should be understood with the mean convex
manifold which comes with the definition. However, being extreme for a curve is the property of the curve alone
(depends only on the ambient manifold).
\end{rmk}

Now, we state the main facts which we use in the following sections.

\begin{lem}\cite{MY2}, \cite{MY3}
Let $M$ be a compact, mean convex $3$-manifold, and $\Gamma\subset\partial M$ be a nullhomotopic simple closed curve.
Then, there exists an area minimizing disk $ D\subset M$ with $\partial  D = \Gamma$. Moreover, all such disks are
properly embedded in $M$ and they are pairwise disjoint. Also, if $\Gamma_1, \Gamma_2 \subset \partial M$ are disjoint
simple closed curves, then the area minimizing disks $ D_1,  D_2$ spanning $\Gamma_1, \Gamma_2$ are also disjoint.
\end{lem}

There is an analogous fact for area minimizing surfaces, too.

\begin{lem} \cite{Fe}, \cite{HSi}, \cite{Wh2}
Let $M$ be a compact, mean convex $3$-manifold, and $\Gamma\subset\partial M$ be a nullhomologous simple closed curve.
Then, there exists a smoothly embedded absolutely area minimizing surface $\Sigma\subset M$ with $\partial \Sigma =
\Gamma$.
\end{lem}

Now, we state a lemma about the limit of area minimizing disks in a mean convex manifold. Note that
we mean that the boundary of the disk is in the boundary of the manifold by being {\em properly
embedded}.

\begin{lem}\cite{HS}
Let $M$ be a compact, mean convex $3$-manifold and let $\{ D_i\}$ be a sequence of properly embedded area minimizing
disks in $M$. Then there is a subsequence $\{ D_{i_j}\}$ of $\{ D_i\}$ such that $ D_{i_j} \rightarrow  \widehat{D}$, a
countable collection of properly embedded area minimizing disks in $\Omega$.
\end{lem}

\noindent {\bf Convention:} Throughout the paper, all the manifolds will be assumed to be compact, orientable, mean
convex and having trivial second homology, i.e. $H_2(M,\BZ)=0$. We will also assume that all the surfaces are
orientable as well.

\section{Generic Uniqueness for Area Minimizing Disks}

In this section, we will prove the generic uniqueness of area minimizing disks for $H$-extreme curves. For this, we
first show that for any nullhomotopic simple closed curve in the boundary of a mean convex $3$-manifold, either there
exists a unique area minimizing disk spanning the curve, or there are two canonical extremal area minimizing disks
which bounds a region containing all other area minimizing disks with same boundary. Similar results also appears in
\cite{MY3}, \cite{Li}, \cite{Wh3} and \cite{Co}.

\begin{lem}
Let $M$ be a compact, orientable, mean convex $3$-manifold with $H_2(M,\BZ)=0$. Let $\Gamma$ be a nullhomotopic (in
$M$) simple closed curve in $\partial M$. Then either there is a unique area minimizing disk $D$ in $M$ with $\partial
D = \Gamma$, or there are two canonical area minimizing disks $ D^+$ and $ D^-$ in $M$ with $\partial  D^\pm = \Gamma$,
and any other area minimizing disk in $M$ with boundary $\Gamma$ must belong to the canonical region $N$ bounded by
$D^+$ and $D^-$ in $M$.
\end{lem}

\begin{pf}
Let $M$ be a mean convex $3$-manifold and let $\Gamma\subset\partial M$ be a nullhomotopic simple closed curve. Take a
small neighborhood $A$ of $\Gamma$ in $\partial M$, which will be a thin annulus where $\Gamma$ is the core. $\Gamma$
separates the annulus $A$ into two parts, say $A^+$ and $A^-$ by giving a local orientation. Define a sequence of
pairwise disjoint simple closed curves $\{\Gamma_i^+\} \subset A^+ \subset \partial M$ such that $\lim\Gamma_i^+ =
\Gamma$. Now, by Lemma 2.1, for any curve $\Gamma_i^+$, there exist an embedded area minimizing disk $ D_i^+$ with
$\partial  D_i^+ = \Gamma_i^+$. This defines a sequence of area minimizing disks $\{ D_i^+\}$ in $M$. By Lemma 2.3,
there exists a subsequence $\{ D_{i_j}^+\}$ converging to a countable collection of area minimizing disks $\widehat{
D}^+$ with $\partial \widehat{ D}^+ = \Gamma$.

We claim that this collection $\widehat{ D}^+$ consists of only one area minimizing disk. Assume
that there are two disks in the collection say $ D_a^+$ and $ D_b^+$, and say $ D_a^+$ is {\em
above} $D_b^+$ (in the positive side of $D_b^+$ in the local orientation). By Lemma 2.1, $D_a^+$
and $D_b^+$ are embedded and disjoint. They have the same boundary $\Gamma\subset \partial M$. $
D_b^+$ is also limit of the sequence $\{ D_i^+\}$. But, since for any area minimizing disk $
D_i^+\subset M$, $\partial  D_i^+ =\Gamma_i^+$ is disjoint from $\partial D_a^+ = \Gamma$, $ D_i^+$
disjoint from $ D_a^+$, again by Lemma 2.1. This means $ D_a^+$ is a barrier between the sequence
$\{ D_i^+\}$ and $ D_b^+$, and so, $ D_b^+$ cannot be limit of this sequence. This is a
contradiction. So $\widehat{ D}^+$ is just one area minimizing disk, say $ D^+$. Similarly,
$\widehat{ D}^- = D^-$.

Now, we claim these area minimizing disks $ D^+$ and $ D^-$ are canonical, depending only on
$\Gamma$ and $M$, and independent of the choice of the sequence $\{\Gamma_i\}$ and $\{ D_i\}$. Let
$\{\gamma_i^+\}$ be another sequence of simple closed curves in $A^+$. Assume that there exists
another area minimizing disk $E^+$ with $\partial E^+ = \Gamma$ and $E^+$ is a limit of the
sequence of area minimizing disks $E_i^+$ with $\partial E_i^+ = \gamma_i^+ \subset A^+$. By Lemma
2.1, $ D^+$ and $E^+$ are disjoint. Then one of them is above the other one. If $ D^+$ is above
$E^+$, then $ D^+$ between the sequence $E_i^+$ and $E^+$. This is because, all $E_i^+$ are
disjoint and above $E^+$ as $\partial E_i^+ = \gamma_i$ are disjoint and above $\Gamma$. Similarly,
$D^+$ is {\em below} $E_i$ for any $i$ (in the negative side of $E_i$ in the local orientation), as
$\partial  D^+ = \Gamma$ is below the curves $\gamma_i^+\subset A^+$. Now, since $ D^+$ is between
the sequence $\{E_i^+\}$ and its limit $E^+$, and $E^+$ and $ D^+$ are disjoint, $ D^+$ will be a
barrier for the sequence $\{E_i^+\}$, and so they cannot limit on $E^+$. This is a contradiction.
Similarly, $ D^+$ cannot be below $E^+$, so they must be same. Hence, $ D^+$ and $ D^-$ are
canonical area minimizing disks for $\Gamma$.

Now, we will show that any area minimizing disk in $M$ with boundary $\Gamma$ must belong to the canonical region $N$
bounded by $D^+$ and $D^-$ in $M$, i.e. $\partial N \supseteq D^+ \cup D^-$ ($H_2(M,\BZ)=0$). Let $E$ be any area
minimizing disk with boundary $\Gamma$. By Lemma 2.1, $E$ is disjoint from $D^+$ and $D^-$. Hence, if $E$ is not in
$N$, then it must be completely outside of $N$. So, $E$ is either above $D^+$ or below $D^-$. However, $D^+ = \lim
D_i^+$ and $\Gamma_i^+ \rightarrow \Gamma$ from above. Moreover, again by Lemma 2.1, $E$ must be disjoint from $D_i^+$.
Hence, $E$ would be a barrier between the sequence $\{ D_{i_j}^+\}$ and $D^+$ like in previous paragraph. This is a
contradiction. Similarly, same is true for $D^-$. Hence, any area minimizing disk in $M$ with boundary $\Gamma$ must
belong to the canonical region $N$ bounded by $D^+$ and $D^-$ in $M$. This also shows that if $D^+ = D^-$, then there
exists a unique area minimizing disk in $M$ with boundary $\Gamma$.
\end{pf}

\begin{rmk}
The results in \cite{MY3}, \cite{Li}, \cite{Wh3} are similar to this one in some sense. In those papers, the authors
show the "strong uniqueness" property, which says that either an $H$-extreme curve bounds more than one {\em minimal
disk} in the mean convex manifold $M$ or there is a unique {\em minimal surface} bounding the curve which is indeed an
area minimizing disk in $M$. Our result is relatively different than the others. In above lemma, we proved that either
there exists a unique area minimizing disk in $M$ bounding the $H$-extreme curve, or there are two canonical extremal
area minimizing disks in $M$ which bounds a region containing all other area minimizing disks with same boundary.
\end{rmk}

Now, we prove the main result of the paper.

\begin{thm}
Let $M$ be a compact, orientable, mean convex $3$-manifold with $H_2(M,\BZ)=0$. Then for a generic nullhomotopic (in
$M$) simple closed curve $\Gamma$ in $\partial M$, there exists a unique area minimizing disk $D$ in $M$ with $\partial
D = \Gamma$. In other words, let $\mathcal{A}$ be the space of nullhomotopic (in $M$) simple closed curves in $\partial
M$ and let $\mathcal{A}'\subset \mathcal{A}$ be the subspace containing the curves bounding a unique area minimizing
disk in $M$. Then, $\mathcal{A}'$ is generic in $\mathcal{A}$, i.e. $\mathcal{A}'$ is countable intersection of open
dense subsets.
\end{thm}

\begin{pf} We will prove this theorem in 2 steps.

\vspace{0.3cm}

\textbf{Claim 1:} $\mathcal{A}'$ is dense in $\mathcal{A}$ as a subspace of $C^0(S^1,\partial M)$
with the supremum metric.

\vspace{0.3cm}

\begin{pf} Let $\mathcal{A}$ be the space of nullhomotopic simple closed curves in $\partial M$. We parametrize
this space with $C^0$ parametrizations, and use supremum metric, i.e. $\mathcal{A}= \{\alpha\in
C^0(S^1,\partial M)\ | \ \alpha(S^1) \mbox{ is an embedding, and nullhomotopic in } M\}$.

Now, let $\Gamma_0\in \mathcal{A}$ be a nullhomotopic simple closed curve in $\partial M$. Since
$\Gamma_0$ is simple, there exists a small closed neighborhood $N(\Gamma_0)$ of $\Gamma_0$ which is
an annulus in $\partial M$. Let $\Gamma:[-\epsilon,\epsilon]\rightarrow \mathcal{A}$ be a small
path in $\mathcal{A}$ through $\Gamma_0$ such that $\Gamma(t)=\Gamma_t$ and $\{\Gamma_t\}$ foliates
$N(\Gamma)$ with simple closed curves $\Gamma_t$. In other words, $\{\Gamma_t\}$ are pairwise
disjoint simple closed curves, and $N(\Gamma_0)=\bigcup_{t\in [-\epsilon,\epsilon]} \Gamma_t$.

By Lemma 3.1, for any $\Gamma_t$ either there exists a unique area minimizing disk $D_t$ in $M$, or
there is a canonical region $N_t$ in $M$ between the canonical area minimizing disks $D_t^+$ and
$D_t^-$. With abuse of notation, if $\Gamma_t$ bounds a unique area minimizing disk $D_t$ in $M$,
define $N_t=D_t$ as a degenerate canonical neighborhood for $\Gamma_t$. Clearly, degenerate
neighborhood $N_t$ means $\Gamma_t$ bounds a unique area minimizing disk, and nondegenerate
neighborhood $N_s$ means that $\Gamma_s$ bounds more than one area minimizing disk. Note that by
Lemma 3.1 and Lemma 2.1, all canonical neighborhoods in the collection are pairwise disjoint.

Now, let $\widehat{N}$ be the union of these canonical neighborhoods $\{N_t\}$, i.e. $\widehat{N} =
\bigcup_{t\in [-\epsilon,\epsilon]}N_t$. Then, $\partial \widehat{N} \supseteq D_\epsilon^+ \cup
N(\Gamma_0) \cup D_{-\epsilon}^-$. Let $p^+$ be a point in $D_\epsilon^+$ and $p^-$ be a point in
$D_{-\epsilon}^-$. Let $\beta$ be a finite curve from $p^+$ to $p^-$ intersecting transversely all
the canonical neighborhoods in the collection $\widehat{N}$.

Now, for each $t\in[-\epsilon,\epsilon]$, we will assign a real number $s_t\geq 0$. Let $I_t = \beta\cap N_t$, and
$s_t$ be the length of $I_t$. Then, if $N_t$ is degenerate (There exists a unique area minimizing disk $D_t$ in $M$ for
$\Gamma_t$), then $s_t$ would be $0$. If $N_t$ is nondegenerate ($\Gamma_t$ bounds more than one area minimizing disk),
then $s_t > 0$. Also, it is clear that for any $t$, $I_t\subset \beta$ and $I_t\cap I_s=\emptyset$ for any $t\neq s$.
Then, $\sum_{t\in[-\epsilon,\epsilon]} s_t < C$ where $C$ is the length of $\beta$. This means for only countably many
$t\in[-\epsilon,\epsilon]$, $s_t > 0$. So, there are only countably many nondegenerate $N_t$ for
$t\in[-\epsilon,\epsilon]$. Hence, for all other $t$, $N_t$ is degenerate. This means there exist uncountably many
$t\in[-\epsilon,\epsilon]$, where $\Gamma_t$ bounds a unique area minimizing disk. Since $\Gamma_0$ is arbitrary, this
proves $\mathcal{A} '$ is dense in $\mathcal{A}$.

\end{pf}

\textbf{Claim 2:} $\mathcal{A}'$ is generic in $\mathcal{A}$.

\vspace{0.3cm}

\begin{pf} We will prove that $\mathcal{A} '$ is countable intersection of open dense subsets of $\mathcal{A}$.
Then the result will follow by Baire category theorem.

Since the space of continuous maps from circle to boundary of $M$, $C^0(S^1,\partial M)$, is
complete with supremum metric, then the closure of $\mathcal{A}$ in $C^0(S^1,\partial M)$,
$\bar{\mathcal{A}}\subset C^0(S^1,\partial M)$, is also complete.

Now, we will define a sequence of open dense subsets $U^i\subset \mathcal{A}$ such that their
intersection will give us $\mathcal{A} '$. Let $\Gamma\in \mathcal{A}$ be a simple closed curve in
$\partial M$. As in the Claim 1, let $N(\Gamma)\subset \partial M$ be a neighborhood of $\Gamma$ in
$\partial M$, which is an open annulus. Then, define an open neighborhood $U_\Gamma$ of $\Gamma$ in
$\mathcal{A}$, such that $U_\Gamma = \{\alpha \in \mathcal{A} \ | \ \alpha(S^1)\subset N(\Gamma), \
\alpha \mbox{ is homotopic to } \Gamma\}$. Clearly, $\mathcal{A}= \bigcup_{\Gamma\in \mathcal{A}}
U_\Gamma$. Now, define a finite curve $\beta_\Gamma$ as in Claim 1, which intersects transversely
all the area minimizing disks bounding the curves in $U_\Gamma$.

Now, for any $\alpha \in U_\Gamma$, by Lemma 3.1, there exists a canonical region $N_\alpha$ in $M$
(which can be degenerate if $\alpha$ bounds a unique area minimizing disk). Let $I_{\alpha,\Gamma}
= N_\alpha \cap \beta_\Gamma$. Then let $s_{\alpha,\Gamma}$ be the length of $I_{\alpha,\Gamma}$
($s_{\alpha,\Gamma}$ is $0$ if $N_\alpha$ degenerate). Hence, for every element $\alpha$ in
$U_\Gamma$, we assign a real number $s_{\alpha,\Gamma} \geq 0$.

Now, we define the sequence of open dense subsets in $U_\Gamma$. Let $U^i_\Gamma = \{\alpha\in
U_\Gamma \ | \ s_{\alpha,\Gamma} < 1/i \  \}$. We claim that $U^i_\Gamma$ is an open subset of
$U_\Gamma$ and $\mathcal{A}$. Let $\alpha\in U^i_\Gamma$, and let $s_{\alpha,\Gamma} = \lambda <
1/i$. So, the interval $I_{\alpha,\Gamma}\subset \beta_\Gamma$ has length $\lambda$. Let $I '
\subset \beta_\Gamma$ be an interval containing $I_{\alpha,\Gamma}$ in its interior, and has length
less than $1/i$. By the proof of Claim 1, we can find two simple closed curves $\alpha^+, \alpha^-
\in U_\Gamma$ with the following properties.

\begin{itemize}

\item $\alpha^\pm$ are disjoint from $\alpha$,

\item $\alpha^\pm$ are lying in opposite sides of $\alpha$ in $\partial M$,

\item $\alpha^\pm$ bounds a unique area minimizing disk $D_{\alpha^\pm}$,

\item $D_{\alpha^\pm} \cap \beta_\Gamma \subset I '$.

\end{itemize}

The existence of such curves is clear from the proof of Claim 1, as if one takes any foliation
$\{\alpha_t\}$ of a small neighborhood of $\alpha$ in $\partial M$, there are uncountably many
curves in the family bounding a unique area minimizing disk, and one can choose sufficiently close
pair of curves to $\alpha$, to ensure the conditions above.

After finding $\alpha^\pm$, consider the open annulus $F_\alpha$ in $\partial M$ bounded by
$\alpha^+$ and $\alpha^-$. Let $V_\alpha = \{ \gamma\in U_\Gamma \ | \ \gamma(S^1)\subset F_\alpha
, \ \gamma \mbox{ is homotopic to } \alpha \}$. Clearly, $V_\alpha$ is an open subset of
$U_\Gamma$. If we can show $V_\alpha\subset U^i_\Gamma$, then this proves $U^i_\Gamma$ is open for
any $i$ and any $\Gamma\in \mathcal{A}$.

Let $\gamma\in V_\alpha$ be any curve, and $N_\gamma$ be its canonical neighborhood given by Lemma
3.1. Since $\gamma(S^1)\subset F_\alpha$, $\alpha^+$ and $\alpha^-$ lie in opposite sides of
$\gamma$ in $\partial M$. This means $D_{\alpha^+}$  and $D_{\alpha^-}$ lie in opposite sides of
$N_\gamma$. By choice of $\alpha^\pm$, this implies $N_\gamma \cap \beta_\Gamma= I_{\gamma,\Gamma}
\subset I '$. So, the length $s_{\gamma,\Gamma}$ is less than $1/i$. This implies $\gamma\in
U^i_\Gamma$, and so $V_\alpha\subset U^i_\Gamma$. Hence, $U^i_\Gamma$ is open in $U_\Gamma$ and
$\mathcal{A}$.

Now, we can define the sequence of open dense subsets. Let $U^i = \bigcup_{\Gamma\in \mathcal{A}} U^i_\Gamma$ be an
open subset of $\mathcal{A}$. Since, the elements in $\mathcal{A} '$ represent the curves bounding a unique area
minimizing disk, for any $\alpha\in \mathcal{A} '$, and for any $\Gamma\in \mathcal{A}$, $s_{\alpha,\Gamma} = 0$. This
means $\mathcal{A}'\subset U^i$ for any $i$. By Claim 1, $U^i$ is open dense in $\mathcal{A}$ for any $i>0$.

As we mention at the beginning of the proof, since the space of continuous maps from circle to
boundary of $M$, $C^0(S^1,\partial M)$ is complete with supremum metric, then the closure
$\bar{\mathcal{A}}$ of $\mathcal{A}$ in $C^0(S^1,\partial M)$ is also complete metric space. Since
$\mathcal{A}'$ is dense in $\mathcal{A}$, it is also dense in $\bar{\mathcal{A}}$. As $\mathcal{A}$
is open in $C^0(S^1,\partial M)$, this implies $U^i$ is a sequence of open dense subsets of
$\bar{\mathcal{A}}$. On the other hand, since $U_1 \supseteq U_2 \supseteq ... \supseteq U_n
\supseteq ...$ and $ \bigcap_{i=1}^\infty U_i = \mathcal{A}'$, $\mathcal{A}'$ is generic in
$\mathcal{A}$.
\end{pf}
\end{pf}

\begin{rmk} Notice that we use the homology condition just to make sure that $D^+\cup D^-$ is a separating sphere in $M$, and hence to define
the canonical region between them in Lemma 3.1. So, if we replace $H_2(M,\BZ)=0$ condition with irreducibility of
$3$-manifold (any embedded $2$-sphere bounds a $3$-ball in $M$), the same proof for Lemma 3.1 and Theorem 3.2 would go
through. In other words, Theorem 3.2 is also true for compact, irreducible, orientable, mean convex $3$-manifolds.
\end{rmk}

\section{Generic Uniqueness for Area Minimizing Surfaces}

In this section, we will prove the generic uniqueness result for $H$-extreme curves in the absolutely area minimizing
case. The technique is basically same with area minimizing disk case. First, we will prove an analogous version of
Lemma 2.1 [MY2, Theorem 6] for absolutely area minimizing surfaces. However, the analogous version of Lemma 2.1 is not
true in general for global version. Hence, we will prove it for a local version which suffices for our purposes. See
Remark 4.1.

\begin{lem} Let $M$ be a compact, orientable, mean convex $3$-manifold with $H_2(M,\BZ)=0$. Let $A$ be an annulus in $\partial M$ whose core
$\Gamma$ is nullhomologous in $M$. If $\Gamma_1$ and $\Gamma_2$ are two disjoint simple closed curves in $A$ which are
homotopic to $\Gamma$ in $A$, then any absolutely area minimizing surfaces $\Sigma_1$ and $\Sigma_2$ in $M$ with
$\partial \Sigma_i = \Gamma_i$ are disjoint, too. Moreover, if $\Sigma$ and $\Sigma '$ are two absolutely area
minimizing surfaces in $M$ where $\partial \Sigma = \partial \Sigma ' = \Gamma$, then they must be disjoint, too.
\end{lem}

\begin{pf}
Let $M$ be a mean convex $3$-manifold, and $A$ is an annulus in $\partial M$ whose core $\Gamma$ is
nullhomologous in $M$. Let $\Gamma_1$ and $\Gamma_2$ are two disjoint simple closed curves in $A$
which are homotopic to $\Gamma$ in $A$. Let $\Sigma_1$ and $\Sigma_2$ be absolutely area minimizing
surfaces in $M$ with $\partial \Sigma_i = \Gamma_i$. We want to show that $\Sigma_1$ and $\Sigma_2$
are disjoint.

Assume on the contrary that $\Sigma_1\cap\Sigma_2 \neq \emptyset$. Now, let $\widehat{N}$ be the
convex hull of $A$ in $M$. Then, by maximum principle, $\Sigma_1$ and $\Sigma_2$ are in
$\widehat{N}$. Moreover, as $\Gamma_1$ separates the annulus $A$, then $\Sigma_1$ is separating in
$\widehat{N}$. Similarly, $\Sigma_2$ is separating, too. Now, if $\Sigma_1\cap\Sigma_2 = \gamma$
where $\gamma$ is a collection of closed curves, then $\Sigma_1$ separates $\Sigma_2$ into two
subsurfaces $S^1_1$ and $S^1_2$ where $\partial S^1_1 = \gamma$ and $\partial S^1_2 = \gamma\cup
\Gamma_1$. Similarly, $\Sigma_2$ separates $\Sigma_1$ into two subsurfaces $S^2_1$ and $S^2_2$
where $\partial S^2_1 = \gamma$ and $\partial S^2_2 = \gamma\cup \Gamma_2$. Now, we will use the
Meeks-Yau exchange roundoff trick to get a contradiction \cite{MY2}.

As $\Sigma_1$ and $\Sigma_2$ are absolutely area minimizing surfaces in $M$, $|S^1_1| = |S^2_1|$
where $ | S |$ is the area of $S$. Now define a new surface by swaping the subsurfaces $S^1_1$ and
$S^2_1$. In other words, let $T_1 = (\Sigma_1 - S^1_1) \cup S^2_1$. As $T_1$ and $\Sigma_1$ have
same area, then $T_1$ is also absolutely area minimizing surface. However, $\gamma$ is a folding
curve in $T_1$ as in \cite{MY2}. This is a contradiction (One can also argue with the regularity of
the absolutely area minimizing surfaces \cite{Fe}). Hence, this shows that $\Sigma_1$ and
$\Sigma_2$ in $M$ with $\partial \Sigma_i = \Gamma_i$ are disjoint absolutely area minimizing
surfaces in $M$.

Now, we will consider same boundary case. Let $A$ and $\Gamma$ be as in the statement of the
theorem. Let $\Sigma$ and $\Sigma '$ be two absolutely area minimizing surfaces where $\partial
\Sigma = \partial \Sigma ' = \Gamma$. Let $\widehat{N}$ be as above. Then, $\Sigma_1$ and
$\Sigma_2$ are separating in $\widehat{N}$. As in the previous paragraph, $\Sigma_1$ and $\Sigma_2$
separates each other, and by swaping argument again, we get a contradiction. The proof follows.
\end{pf}

\begin{rmk}
The techniques for Lemma 2.1 (or [MY2, Theorem 6]) is not working for an analogous theorem in absolutely area
minimizing surfaces case in general. In other words, if we just assume $\Gamma_1\cap\Gamma_2=\emptyset$, and not
require them to be in the annulus $A$, then the techniques of the above lemma do not apply. This is because if, for
example, $\Gamma_1$ or $\Gamma_2$ are not separating in $\partial M$, then the intersection of absolutely area
minimizing surfaces $\Sigma_1$ and $\Sigma_2$ might contain a nonseparating curve $\gamma$ in one of the surfaces, say
$\Sigma_1$. Hence, we cannot make any surgery there because $\gamma$ may not bound a subsurface in $\Sigma_1$. So, we
went to a local version of this theorem (which is enough for our purposes) by restricting $\partial M$ to a small
subannnulus $A$ in $\partial M$ to make sure that each essential curve is separating in $A$, and we can make surgery in
the intersection of surfaces.
\end{rmk}

Now, we will give a generalization of Lemma 3.1 in absolutely area minimizing surface case.

\begin{lem}
Let $M$ be a compact, orientable, mean convex $3$-manifold with $H_2(M,\BZ)=0$. Let $\Gamma$ be a nullhomologous (in
$M$) simple closed curve in $\partial M$. Then either there is a unique absolutely area minimizing surface $\Sigma$ in
$M$ with $\partial \Sigma = \Gamma$, or there are uniquely defined two canonical extremal absolutely area minimizing
surfaces $\Sigma^+$ and $\Sigma^-$ in $M$ with $\partial \Sigma^\pm = \Gamma$, and any other absolutely area minimizing
surface in $M$ with boundary $\Gamma$ must belong to the canonical region $N$ bounded by $\Sigma^+$ and $\Sigma^-$ in
$M$.
\end{lem}

\begin{pf}
Let $M$ be a mean convex $3$-manifold and let $\Gamma\subset\partial M$ be a nullhomologous simple closed curve. Take a
small neighborhood $A$ of $\Gamma$ in $\partial M$, which will be a thin annulus where $\Gamma$ is the core. $\Gamma$
separates the annulus $A$ into two parts, say $A^+$ and $A^-$ by giving a local orientation. Define a sequence of
pairwise disjoint simple closed curves $\{\Gamma_i^+\} \subset A^+ \subset \partial M$ such that $\lim\Gamma_i^+ =
\Gamma$. Now, by Lemma 2.2, for any curve $\Gamma_i^+$, there exist an embedded absolutely area minimizing surface
$\Sigma_i^+$ with $\partial \Sigma_i^+ = \Gamma_i^+$. This defines a sequence of absolutely area minimizing surfaces
$\{ \Sigma_i^+\}$ in $M$. By \cite{Fe}, there exists a subsequence $\{\Sigma_{i_j}^+\}$ converging to an absolutely
area minimizing surface $\Sigma^+$ with $\partial \Sigma^+ = \Gamma$. Similarly, by defining a similar sequence
$\{\Gamma_i^-\}$ in $A^-$ and similar construction, an absolutely area minimizing surface $\Sigma^-$ with $\partial
\Sigma^- = \Gamma$ can be defined.

Now, we claim these absolutely area minimizing surfaces $\Sigma^+$ and $\Sigma^-$ are canonical, depending only on
$\Gamma$ and $M$, and independent of the choice of the sequence $\{\Gamma_i\}$ and $\{ \Sigma_i\}$. Let
$\{\gamma_i^+\}$ be another sequence of simple closed curves in $A^+$. Assume that there exists another absolutely area
minimizing surface $S^+$ with $\partial S^+ = \Gamma$ and $S^+$ is a limit of the sequence of absolutely area
minimizing surfaces $S_i^+$ with $\partial S_i^+ = \gamma_i^+ \subset A^+$. As $\Sigma^+$ and $S^+$ are absolutely area
minimizing surfaces with same boundary $\Gamma$, they are disjoint by Lemma 4.1. Then one of them is {\em above} the
other one (in the positive side of the other one in the local orientation). If $\Sigma^+$ is above $S^+$, then
$\Sigma^+$ between the sequence $S_i^+$ and $S^+$. This is because, all $S_i^+$ are disjoint and above $S^+$ as
$\partial S_i^+ = \gamma_i$ are disjoint and above $\Gamma$. Similarly, $\Sigma^+$ is below $S_i$ for any $i$, as
$\partial \Sigma^+ = \Gamma$ is below the curves $\gamma_i^+\subset A^+$. Now, since $\Sigma^+$ is between the sequence
$\{S_{i_j}^+\}$ and its limit $S^+$, and $S^+$ and $\Sigma^+$ are disjoint, $\Sigma^+$ will be a barrier for the
sequence $\{S_{i_j}^+\}$, and so they cannot limit on $S^+$. This is a contradiction. Similarly, $\Sigma^+$ cannot be
below $S^+$, so they must be same. Hence, $\Sigma^+$ and $\Sigma^-$ are canonical absolutely area minimizing surfaces
for $\Gamma$.

Now, we will show that any absolutely area minimizing surface in $M$ with boundary $\Gamma$ must belong to the
canonical region $N$ bounded by $\Sigma^+$ and $\Sigma^-$ in $M$, i.e. $\partial N \supseteq \Sigma^+ \cup\Sigma^-$
($H_2(M,\BZ)=0$). Let $T$ be any absolutely area minimizing surface with boundary $\Gamma$. By Lemma 4.1, $T$ is
disjoint from $\Sigma^+$ and $\Sigma^-$. Hence, if $T$ is not in $N$, then it must be completely outside of $N$. So,
$T$ is either above $\Sigma^+$ or below $\Sigma^-$. Assume $T$ is above $\Sigma^+$. However, $\Sigma^+ = \lim
\Sigma_{i_j}^+$ and $\Gamma_{i_j}^+ \rightarrow \Gamma$ from above. Moreover, again by Lemma 4.1, $T$ must be disjoint
from $\Sigma_{i_j}^+$. Hence, $T$ is a barrier between the subsequence $\Sigma_{i_j}^+$ and its limit $\Sigma$. Like in
previous paragraph, this is a contradiction. Similarly, the same is true for $\Sigma^-$. Hence, any absolutely area
minimizing surface in $M$ with boundary $\Gamma$ must belong to the canonical region $N$ bounded by $\Sigma^+$ and
$\Sigma^-$ in $M$. This also shows that if $\Sigma^+ =\Sigma^-$, then there exists a unique absolutely area minimizing
surface in $M$ with boundary $\Gamma$.
\end{pf}

Now, we can prove the generic uniqueness result for absolutely area minimizing surfaces.

\begin{thm}
Let $M$ be a compact, orientable, mean convex $3$-manifold with $H_2(M,\BZ)=0$. Then for a generic nullhomologous (in
$M$) simple closed curve $\Gamma$ in $\partial M$, there exists a unique absolutely area minimizing surface $\Sigma$ in
$M$ with $\partial \Sigma = \Gamma$. In other words, let $\mathcal{A}$ be the space of nullhomologous (in $M$) simple
closed curves in $\partial M$ and let $\mathcal{A}'\subset \mathcal{A}$ be the subspace containing the curves bounding
a unique absolutely area minimizing surface in $M$. Then, $\mathcal{A}'$ is generic in $\mathcal{A}$, i.e.
$\mathcal{A}'$ is countable intersection of open dense subsets.
\end{thm}

\begin{pf}
The idea is basically same with Theorem 3.2. We will imitate the same proof in this context. Again, we will prove this
theorem in 2 steps.

\vspace{0.3cm}

\textbf{Claim 1:} $\mathcal{A}'$ is dense in $\mathcal{A}$ as a subspace of $C^0(S^1,\partial M)$
with the supremum metric.

\vspace{0.3cm}

\begin{pf} Let $\mathcal{A}$ be the space of nullhomologous simple closed curves in $\partial M$. We parametrize
this space with $C^0$ parametrizations, and use supremum metric, i.e. $\mathcal{A}= \{\alpha\in
C^0(S^1,\partial M)\ | \ \alpha(S^1) \mbox{ is an embedding, and nullhomologous in } M\}$.

Now, let $\Gamma_0\in \mathcal{A}$ be a nullhomologous simple closed curve in $\partial M$. As in
the proof of Theorem 3.2, let $N(\Gamma_0)$ be an annulus in $\partial M$ and Let
$\Gamma:[-\epsilon,\epsilon]\rightarrow \mathcal{A}$ foliates $N(\Gamma)$ with simple closed curves
$\Gamma_t$.

By Lemma 4.2, for any $\Gamma_t$ either there exists a unique absolutely area minimizing surface
$\Sigma_t$ in $M$, or there is a canonical region $N_t$ in $M$ between the canonical area
minimizing surfaces $\Sigma_t^+$ and $\Sigma_t^-$. As in the proof of Theorem 3.2, if $\Gamma_t$
bounds a unique absolutely area minimizing surface $\Sigma_t$ in $M$, define $N_t=\Sigma_t$ as a
degenerate canonical neighborhood for $\Gamma_t$. Clearly, degenerate neighborhood $N_t$ means
$\Gamma_t$ bounds a unique absolutely area minimizing surface, and nondegenerate neighborhood $N_s$
means that $\Gamma_s$ bounds more than one absolutely area minimizing surface. Note that by Lemma
4.1, all canonical neighborhoods in the collection are pairwise disjoint.

Like before, let $\widehat{N} = \bigcup_{t\in [-\epsilon,\epsilon]} N_t$.  Let $p^+$ be a point in
$\Sigma_\epsilon^+$ and $p^-$ be a point in $\Sigma_{-\epsilon}^-$. Let $\beta$ be a finite curve
from $p^+$ to $p^-$ intersecting transversely all the canonical neighborhoods in the collection
$\widehat{N}$. For each $t\in[-\epsilon,\epsilon]$, assign a real number $s_t$ to be the length of
$I_t = \beta \cap N_t$. Clearly if $N_t$ is nondegenerate ($\Gamma_t$ bounds more than one
absolutely area minimizing surface), then $s_t > 0$. Then, $\sum_{t\in[-\epsilon,\epsilon]} s_t <
C$ where $C$ is the length of $\beta$. This means for only countably many
$t\in[-\epsilon,\epsilon]$, $s_t > 0$. So, there are only countably many nondegenerate $N_t$ for
$t\in[-\epsilon,\epsilon]$. Hence, for all other $t$, $N_t$ is degenerate. This means there exist
uncountably many $t\in[-\epsilon,\epsilon]$, where $\Gamma_t$ bounds a unique absolutely area
minimizing surface. Since $\Gamma_0$ is arbitrary, this proves $\mathcal{A} '$ is dense in
$\mathcal{A}$.
\end{pf}

\textbf{Claim 2:} $\mathcal{A}'$ is generic in $\mathcal{A}$.

\vspace{0.3cm}

\begin{pf}
Let $\mathcal{A}$ be as in the proof of Theorem 3.2. Again, we will define a sequence of open dense
subsets $U^i\subset \mathcal{A}$ such that their intersection will give us $\mathcal{A} '$. Let
$\Gamma\in \mathcal{A}$ be a simple closed curve in $\partial M$. As in the Claim 1, let
$N(\Gamma)\subset
\partial M$ be a neighborhood of $\Gamma$ in $\partial M$, which is an open annulus. Then, define
an open neighborhood $U_\Gamma$ of $\Gamma$ in $\mathcal{A}$, such that $U_\Gamma = \{\alpha \in
\mathcal{A} \ | \ \alpha(S^1)\subset N(\Gamma), \ \alpha \mbox{ is homotopic to } \Gamma\}$.
Clearly, $\mathcal{A}= \bigcup_{\Gamma\in \mathcal{A}} U_\Gamma$.  Now, define a finite curve
$\beta_\Gamma$ as in Claim 1, which intersects all the absolutely area minimizing surfaces bounding
the curves in $U_\Gamma$.

Now, for any $\alpha \in U_\Gamma$, by Lemma 4.2, there exists a canonical region $N_\alpha$ in
$M$. Let $I_{\alpha,\Gamma} = N_\alpha \cap \beta_\Gamma$. Then let $s_{\alpha,\Gamma}$ be the
length of $I_{\alpha,\Gamma}$ ($s_{\alpha,\Gamma}$ is $0$ if $N_\alpha$ degenerate). Now, we define
the sequence of open dense subsets in $U_\Gamma$.

Let $U^i_\Gamma = \{\alpha\in U_\Gamma \ | \ s_{\alpha,\Gamma} < 1/i \  \}$. We claim that
$U^i_\Gamma$ is an open subset of $U_\Gamma$ and $\mathcal{A}$. Let $\alpha\in U^i_\Gamma$, and let
$s_{\alpha,\Gamma} = \lambda < 1/i$. So, the interval $I_{\alpha,\Gamma}\subset \beta_\Gamma$ has
length $\lambda$. Let $I ' \subset \beta_\Gamma$ be an interval containing $I_{\alpha,\Gamma}$ in
its interior, and has length less than $1/i$. Now, let $\alpha^+$ and $\alpha^-$ be as in the proof
of Theorem 3.2. Consider the open annulus $F_\alpha$ in $\partial M$ bounded by $\alpha^+$ and
$\alpha^-$. Let $V_\alpha = \{ \gamma\in U_\Gamma \ | \ \gamma(S^1)\subset F_\alpha , \ \gamma
\mbox{ is homotopic to } \alpha \}$. Clearly, $V_\alpha$ is an open subset of $U_\Gamma$. If we can
show $V_\alpha\subset U^i_\Gamma$, then this proves $U^i_\Gamma$ is open for any $i$ and any
$\Gamma\in \mathcal{A}$.

Let $\gamma\in V_\alpha$ be any curve, and $N_\gamma$ be its canonical neighborhood given by Lemma
4.2. Since $\gamma(S^1)\subset F_\alpha$, $\alpha^+$ and $\alpha^-$ lie in opposite sides of
$\gamma$ in $\partial M$. This means $\Sigma_{\alpha^+}$  and $\Sigma_{\alpha^-}$ lie in opposite
sides of $N_\gamma$. By choice of $\alpha^\pm$, this implies $N_\gamma \cap \beta_\Gamma=
I_{\gamma,\Gamma} \subset I '$. So, the length $s_{\gamma,\Gamma}$ is less than $1/i$. This implies
$\gamma\in U^i_\Gamma$, and so $V_\alpha\subset U^i_\Gamma$. Hence, $U^i_\Gamma$ is open in
$U_\Gamma$ and $\mathcal{A}$. The remaining part of the proof is just like Theorem 3.2.
\end{pf}

\end{pf}

\section{Concluding Remarks}

In this paper, we showed that for a generic nullhomotopic, simple closed curve in the boundary of a mean convex
$3$-manifold $M$, there exists a unique area minimizing disk in $M$. We also prove a similar theorem for absolutely
area minimizing surfaces. In many sense, the techniques used in this paper are purely topological and simple. They are
quite original and can be applied to many similar settings.

Note that all the results of this paper for compact $3$-manifolds with mean convex boundary. For the noncompact case,
like in \cite{MY3} and \cite{HS}, with additional condition of being homogeneously regular on the manifold $M$, all the
results of this paper will go through easily by using the analogous theorems from the same references.

There have been many embeddedness and uniqueness results for the Plateau problem. In extreme and $H$-extreme curve
case, there have been many embeddedness results like \cite{TT}, \cite{AS}, \cite{MY1}. There are also "strong
uniqueness" results for $H$-extreme curves like \cite{MY3}, \cite{Li}, \cite{Wh3}. However, those results do not say
anything about the number of area minimizing disks bounded by an $H$-extreme curve. In those papers, authors gave a
dichotomy that either an $H$-extreme curve bounds more than one minimal disk, or the only minimal surface bounded by
that curve is an area minimizing disk. One should not combine this result with ours in a wrong way. Our result tells
that a generic $H$-extreme curve bounds a unique area minimizing disk. However, bounding a unique area minimizing disk
does not prohibit to bound other minimal surfaces. So, it is not true that for a generic $H$-extreme curve, the only
minimal surface bounded by that curve is an area minimizing disk.

On the other hand, generic uniqueness for area minimizing disks and generic uniqueness for absolutely area minimizing
surfaces might sound contradicting at the first glance. This is because if we have an absolutely area minimizing
surface (which is not a disk) in $M$, we can construct two different area minimizing disks in different sides of the
surface. There are two points to consider here. The first obvious thing is that an absolutely area minimizing surface
might be a disk. The other less obvious fact is that having two different disks in different sides of the surface does
not mean that the curve has more than one area minimizing disk. This is because they are area minimizing disks in that
part of $M$, not the whole $M$. The area minimizing disk in $M$ could be completely different than the others, and it
still can be a unique area minimizing disk in $M$ bounded by that curve.

Another important point here is that these techniques may not work for surfaces which are area minimizing in a fixed
topological class. If they are not absolutely area minimizing in the homology class, or area minimizing disk, then
Lemma 2.1 and its local generalization Lemma 4.1 are not true in general. One should keep in mind that two {\em just}
minimal surfaces with same extreme boundary curve can intersect in a certain way, but two area minimizing disks, or two
absolutely area minimizing surfaces must stay disjoint because of area constraints (intersection implies area
reduction). In those lemmas, we are essentially using Meeks-Yau exchange roundoff trick, and a surgery argument.
However, two surfaces which are area minimizing in a fixed topological class may not give a surface in the same
topological class after surgery. Hence, the key point in our technique (disjointness for the summation argument) fails
in this case. However, as we pointed out in the introduction, White \cite{Wh1} already gave a strong generic uniqueness
result for this case in any dimension and codimension with some smoothness condition.


\begin{thebibliography}{MSY}

\bibitem[AS]{AS} F.J. Almgren, L. Simon, {\em Existence of embedded solutions of Plateau's
problem}, Ann. Scuola Norm. Sup. Pisa Cl. Sci. (4) {\bf 6} (1979) no. 3, 447--495.

\bibitem[Al]{Al} H.W. Alt, {\em Verzweigungspunkte von $H$-Flachen. II}, Math. Ann. {\bf 201} (1973), 33--55.

\bibitem[Co]{Co} B. Coskunuzer, {\em Mean Convex Hulls and Least Area Disks spanning Extreme Curves},
Math Z. {\bf 252} (2006) 811--824.

\bibitem[Do]{Do} J. Douglas, {\em Solution of the problem of Plateau}, Trans. Amer. Math. Soc. {\bf 33} (1931)
263-321.

\bibitem[EWW]{EWW} T. Ekholm, B. White, D. Wienholtz, {\em Embeddedness of minimal surfaces with total boundary curvature at most $4\pi$},
Ann. of Math. (2) {\bf 155} (2002), no. 1, 209--234.

\bibitem[Fe]{Fe} H. Federer, {\em Geometric measure theory}, Springer-Verlag, New York 1969.

\bibitem[Gu]{Gu} R.D. Gulliver, {\em Regularity of minimizing surfaces of prescribed mean curvature}, Ann. of
Math. (2) {\bf 97} (1973) 275--305.

\bibitem[GS]{GS} R. Gulliver, J. Spruck, {\em On embedded minimal surfaces}, Ann. of Math. (2) {\bf 103} (1976) 331--347.

\bibitem[HS]{HS} J. Hass, P. Scott, {\em The existence of least area surfaces in $3$-manifolds},
Trans. Amer. Math. Soc. {\bf 310} (1988) no. 1, 87--114.

\bibitem[HSi]{HSi} R. Hardt, L. Simon, {\em Boundary regularity and embedded solutions for the oriented Plateau
problem}, Ann. of Math. (2) {\bf 110} (1979), no. 3, 439--486.

\bibitem[Li]{Li} F.H. Lin, {\em Plateau's problem for $H$-convex curves}, Manuscripta Math. {\bf 58} (1987) 497--511.

\bibitem[M]{M} F. Morgan, {\em Almost every curve in $R\sp{3}$ bounds a unique area minimizing surface}, Invent.
Math. {\bf 45} (1978) no. 3, 253--297.

\bibitem[Mo]{Mo} C.B. Morrey, {\em The problem of Plateau on a Riemannian manifold}, Ann. of Math. (2)
{\bf 49} (1948) 807--851.

\bibitem[MY1]{MY1} W. Meeks and S.T. Yau, {\em Topology of three-dimensional manifolds and the embedding
problems in minimal surface theory}, Ann. of Math. (2) {\bf 112} (1980) 441--484.

\bibitem[MY2]{MY2} W. Meeks and S.T. Yau, {\em The classical Plateau problem and the topology of three
manifolds}, Topology {\bf 21} (1982) 409--442.

\bibitem[MY3]{MY3} W. Meeks and S.T. Yau, {\em The existence of embedded minimal surfaces and the problem
of uniqueness}, Math. Z. {\bf 179} (1982) no. 2, 151--168.

\bibitem[Ni]{Ni} J.C.C. Nitsche, {\em A new uniqueness theorem for minimal surfaces}, Arch. Rational
Mech. Anal. {\bf 52} (1973) 319--329.

\bibitem[Os]{Os} R. Osserman, {\em A proof of the regularity everywhere of the classical solution to Plateau's
problem}, Ann. of Math. (2) {\bf 91} (1970) 550--569.

\bibitem[Ra]{Ra} T. Rado, {\em On Plateau's problem}, Ann. of Math. (2) {\bf 31} (1930) no. 3, 457--469.

\bibitem[TT]{TT} F. Tomi, A.J. Tromba, {\em Extreme curves bound embedded minimal surfaces of the type of the disc},
Math. Z. {\bf 158} (1978) no. 2, 137--145.

\bibitem[Tr]{Tr} A.J. Tromba, {\em The set of curves of uniqueness for Plateau's problem has a dense interior},
Lecture Notes in Math., Vol. 597, pp. 696--706, Springer, Berlin, 1977.

\bibitem[Wh1]{Wh1} B. White,  {\em The space of $m$-dimensional surfaces that are stationary for a parametric elliptic
functional} Indiana Univ. Math. J. {\bf 36} (1987), no. 3, 567--602.

\bibitem[Wh2]{Wh2} B. White, {\em Existence of smooth embedded surfaces of prescribed genus that
minimize parametric even elliptic functionals on $3$-manifolds}, J. Diff. Geom. {\bf 33} (1991) no.
2, 413--443.

\bibitem[Wh3]{Wh3} B. White, {\em On the topological type of minimal submanifolds}, Topology {\bf 31} 445--448 (1992).


\end{thebibliography}
\end{document}